\documentclass[12pt]{article}
\usepackage{amsmath,amscd,amsthm,a4,amssymb,amsfonts}
\usepackage{mathptmx}      % use Times fonts if available on your TeX system
\usepackage{pst-plot}
\usepackage{color}
\usepackage{latexsym}

\newcommand{\ds}{\displaystyle}
{}
{}
{}
\newtheorem{definition}{Definition}
\newtheorem{theorem}{Theorem}
\theoremstyle{remark}
{}
\newtheorem{Remarks}{Remarks}{}

\chardef\No=242

\begin{document}

\begin{center}
{\Large \bf Hardy's type inequality for the over critical exponent associated with the
Dunkl transform}
\end{center}

\begin{center}
\textbf{Rahmouni Atef}\\
\textit{University of Carthage, Faculty of Sciences of
Bizerte}\\
\textit{Department of Mathematics Bizerte  7021 Tunisia.}\\
Atef.Rahmouni@fsb.rnu.tn\\
\end{center}
\begin{abstract}
For the Hardy space $H^p(\mathbb{R}^{d})$, $ 0<p\leq 1,$ we shall
improve  a Hardy's type inequality associated with  Dunkl
transform respect to  the measures $d\mu_{k}$ homogeneous  of degree $\gamma ,$ on the strip $(2\gamma+d)(2-p)\leq\sigma<2\gamma+d+p(N+1),$
where $N = [(2\gamma+d)( 1/ p -1)]$ is the greatest integer not exceeding
$(2\gamma+d)( 1/ p -1).$ 
\end{abstract}

\begin{quote}\small

\end{quote}
\begin{quote}\small
{\it keywords:} Hardy space, Hardy's type inequality, Dunkl transform.
\end{quote}

\vspace{1mm} \noindent

\section{ Introduction}

In recent years the topic of Hardy type inequalities and their applications
seem to have grown more and more popular.
Although Hardy's original result dates back to the 1920's, some new
versions are stated and old ones are still being improved almost a century
later. One of the reasons for the popularity of Hardy type inequalities is their
usefulness in various applications.

The first definition  of Hardy spaces was in terms of analytic functions
in the unit disc and their boundary values. In the last two decades,
the theory was developed in $\mathbb{R}^{d}$ by real variable methods
like Poisson integrals, Riesz transforms, and maximal functions.
The subsequent discovery of the atomic decomposition theory
of $H^{p}(\mathbb{R}^{d})$ spaces marks an important step of
further developments on its real variable theory.
Using the grand maximal function, R. Coifman \cite{CO}
first shows that an element in $H^{p}(\mathbb{R}^{d})$
can decomposed into a sum of a series of basic elements.
Then the study on some analytic problems on $H^{p}(\mathbb{R}^{d})$
is summed up to investigate some properties of these basic elements,
and therefore the problems because quite simple. Taibleson and Weiss
\cite{TW} gave the definition of molecules belonging to $H^p,$ and
showed that every molecule is in $H^p$ with continuous embedding
map. By the atomic decomposition and the molecule
characterization, the proof of $H^p$ boundedness of the operators
on Hardy space becomes easier. The theory of $H^p$ have been
extensively studied in \cite{FS,S3,GR}

In the setting of the euclidian case, Hardy's inequality for
Fourier transform asserts that for all $f\in H^p(\mathbb{R}^d),$
$0<p\leq 1$
\begin{equation} \label{A}
\int_{\mathbb{R}^d}{|\widehat{f}(\xi)|^p\over |\xi|^{d(2-p)}}
d\xi\leq \|f\|^p_{H^p{(\mathbb{R}^d})},\qquad 0<p\leq1
\end{equation}
where $H^p(\mathbb{R}^d)$ indicates the real Hardy space.
Recently, an extension has been given by \cite{AM1}, the latter
establish a Hardy's type inequality associated with the euclidean
Fourier transform for over critical exponent $\sigma>d(2-p).$

In the same context  we prove a Hardy's type inequality associated with Dunkl transform.
We consider the differential-difference operators  $T_{j},~j=1,...,d,$ on $\mathbb{R}^{d}$ introduced by C.F. Dunkl in 
\cite{Dunkl} and called Dunkl operators in the
literature, associated with the finite reflection group $G$ and the multiplicity function $k,$ are given for a function $f$ of class $C^1$ on
$\mathbb{R}^{d}$ by
\begin{eqnarray*}
T_{j}f(y)&=& \frac{\partial}{\partial y_{j}}f(y)+\sum_{\alpha\in\mathcal{R}_{+}}k(\alpha)
\alpha_{j}\frac{f(y)-f(\sigma_{\alpha}y)}{\prec\alpha,y\succ}.
 \end{eqnarray*}
For $y\in \mathbb{R}^{d},$ the initial problem $T_{j}u(., y)(x) = y_{j}u(x, y);~j=1,...,d$
with $u(0, y) = 1$ admits a unique analytic solution on $\mathbb{R}^{d},$ which will be denoted by
$E_{k}(x, y)$  and called Dunkl kernel \cite{Dunkl1}.
This kernel has a unique analytic extension to $\mathbb{C}^{d}\times \mathbb{C}^{d}.$ The Dunkl kernel has the Laplace-type representation \cite{Rosler}
$$E_{k}(x, y)=\int_{\mathbb{R}^{d}}e^{\prec y,z \succ}d\Gamma_{x}(z);~~ x\in \mathbb{R}^{d},~y\in \mathbb{C}^{d},$$
where $\prec y,z \succ:=\sum_{j=1}^{d}y_{j}z_{j}$ and  $\Gamma_{x}$ is a probability measure on $\mathbb{R}^{d},$ such that $supp(\Gamma_{x})\subset \{z\in\mathbb{R}^{d}:~ |z|\leq|x|\}.$\\
In particular cases, we have
\begin{equation}\label{3}
|E_{k}(x, y)|\leq 1,~~x,y\in\mathbb{R}^{d}.
\end{equation}
The Dunkl kernel gives rise to an integral transform, which is called Dunkl transform on $\mathbb{R}^{d},$ and it was introduced by Dunkl \cite{Dunkl2}, where already many basic properties were established. Dunkl's
results were completed and extended later on by de Jeu \cite{Jeu}. The Dunkl transform of a function $f
\in L^{1} (\mathbb{R}^{d})$ is
$$ \mathcal{F}_{\mathcal{D}}(f)(x):=c_{k}\int_{\mathbb{R}^{d}}E_{k}(-ix,y)f(y)d\mu_{k}(y),~~x\in\mathbb{R}^{d},$$
where  $c_{k}$ is the Mehta-type constant given by  $c_{k}=\Big(\int_{\mathbb{R}^{d}}e^{-|y|^{2}/2}d\mu_{k}(y)\Big)^{-1}.$
We denote by $d\mu_{k}$ the measure on $\mathbb{R}^d$ given by $d\mu_{k}(y) := w_{k}(y)dy$ and by $L^p_{k}(\mathbb{R}^{d}),$ $0< p\leq\infty,$ the space of measurable functions $f$ on $\mathbb{R}^d,$ such that
$$\|f\|_{L^{p}_{k}}=\left(\ds\int_{\mathbb{R}^{d}}|f(y)|^pd\mu_{k}(y)\right)^{1\over
p},\,\mbox{if } p>0 \quad   \mbox{and }\quad \|f\|_{L^{\infty}_{k}}=ess \sup_{y\in\mathbb{R}^{d}}|f(y)|.$$

In this paper, we obtain an improved Hardy's type inequality associated with the Dunkl
transform. So, for the Hardy space $H^p(\mathbb{R}^{d}),~0 <p\leq 1$
we establish a Hardy's type inequality for the strip
$(2\gamma+d)(2-p)\leq\sigma<2\gamma+d+p(N+1).$

Throughout this paper, $C$ stands for a positive constant that can
be changed from line to line.

\section{Preliminaries : \small {(Reflection groups, Root systems and Multiplicity functions)}} 

Let us begin to recall some results concerning the root systems. A useful reference for this topic is the book by Humphreys \cite{Humphreys}.

We consider $\mathbb{R}^{d}$ with the euclidean inner product $\prec.,. \succ$ and norm $|y| :=\sqrt{\prec y,y \succ}.$\\
For $\alpha\in \mathbb{R}^{d}\setminus \{0\},$ let $\sigma_{\alpha}$ be the reflection in the hyperplane
$H_{\alpha}\subset \mathbb{R}^{d}$ orthogonal to $\alpha,$ i.e
$$\sigma_{\alpha} y:=y-\frac{2\prec \alpha,y \succ}{|\alpha|^{2}}\alpha.$$
A finite set $\mathcal{R} \subset\mathbb{R}^{d}\setminus \{0\}$ is called a root system, if
$\mathcal{R} \cap \mathbb{R}.\alpha=\{-\alpha,\alpha\}$ and $\sigma_{\alpha}\mathcal{R}=\mathcal{R}$
for all $\alpha\in\mathcal{R}.$ We assume that it is normalized by $|\alpha|^{2} = 2$ for all $\alpha\in\mathcal{R}.$ For a root system $\mathcal{R},$ the
reflections $\sigma_{\alpha},\alpha\in \mathcal{R},$ generate a finite group $G\subset O(d),$ the reflection group associated with $\mathcal{R}.$
All reflections in $G$ correspond to suitable pairs of roots. For a given $\beta\in \mathbb{R}^{d}\setminus
\bigcup_{\alpha\in\mathcal{R}}H_{\alpha},$ we
fix the positive subsystem $\mathcal{R}_{+}:= \{\alpha\in\mathcal{R}~ :\prec \alpha,\beta\succ > 0\}.$
Then, for each $\alpha \in \mathcal{R}$ either $\alpha \in \mathcal{R}_{+}$ or $\alpha \in \mathcal{R}_{-}.$
Let $k :\mathcal{R}\rightarrow \mathbb{C}$ be a multiplicity function on $\mathcal{R}$ (i.e. a function which is constant on the orbits under the action of $G$). For abbreviation, we introduce the index
$$\gamma=\gamma_{k}:=\sum_{\alpha\in\mathcal{R}_{+}}k(\alpha) .$$
Throughout the paper, we assume that the multiplicity is non-negative, that is, $k(\alpha)\geq 0$ for all
$\alpha \in \mathcal{R}.$  Moreover, let $w_{k}$ denote the weight function
$$w_{k}(y):=\prod_{\alpha\in\mathcal{R}_{+}}|\prec\alpha,y\succ|^{2k(\alpha)},~~ y\in\mathbb{R}^{d},$$
which is $G$-invariant and homogeneous of degree $2\gamma.$

\section{Hardy-type inequality}

The atom decomposition theory of $H^{p}(\mathbb{R}^{d})$ spaces
marks an important step of further developments on its real variable
theory. Using the grand maximal function, R. R. Coifman \cite{CO}
first shows that an element in $H^{p}(\mathbb{R})$ can decomposed
into a sum of a series of basic elements. Then the study on some
analytic problems on $H^{p}(\mathbb{R}^{d})$ is summed up to
investigate some properties of these basic elements, and therefore
the problems because quite simple. These basic elements are called
atoms. Let us now make the definition of an atom.
\begin{definition}
Let $0 < p \leq 1 \leq q \leq \infty $  with $p\neq q.$ A function
$a(x)\in L^{q}(\mathbb{R}^{d})$ is called a $(p, q, s)$-atom with the center at
$x_{0},$ if it satisfies the following conditions
\begin{itemize}
\item[(i)] Supp $ a \subset B(x_{0},r)$
\item[(ii)]$ \|a\|_{L^{q}(\mathbb{R}^{d})}\leq |B(x_{0},r)|^{\frac{1}{q}-\frac{1}{p}};$
\item[(iii)]$ \int_{\mathbb{R}^{d}} a(y)y^{\ell}d\mu_{k}(y)=0, $
for all monomials $y^{\ell}$ with $|\ell|\leq s $ with $s \geq
N=\Big[(2\gamma+d)(\frac{1}{p}-1)\Big],$ where [ . ] denotes, as usual, the
``greatest integer not exceeding" function.
\end{itemize}
\end{definition}
Here, $(i)$ means that an atom must be a function with compact
support, $(ii)$ is the size condition of atoms, and $(iii)$ is
called the cancelation moment condition. Moreover, $B(x_{0},r)$ is
the ball centered at $x_{0}$ with radius $r.$ Clearly, $a(p, \infty,s)$
atom must be $a(p, q, s)$ atom, if $p<q<\infty.$

Using the atomic decomposition, we define the Hardy space
$H^{p}(\mathbb{R}^{d})$ to be the collection of functions $f$
satisfying $f = \sum_{j =0}^{\infty}\beta_{j}a_{j},$ where $a_{j}$
are $H^{p}(\mathbb{R}^{d})$-atoms and $\beta_{j}$ is a sequence of
complex numbers with $\sum_{j=0}^{\infty}|\beta_{j}|^{p}< \infty.$ $
H^{p}(\mathbb{R}^{d})$ is equipped with a norm as follows
$$\|f\|_{H^{p}(\mathbb{R}^{d})}=\inf\Big\{\sum_{j=0}^{\infty}
|\beta_{j}|^{p}\Big\},$$ where the infimum is taken over all atoms
decompositions of $f.$

\begin{theorem}\label{t3.3.11}
Let $0<p \leq 1,$ and  $ N=[(2\gamma+d)({1/p}-1)], $ the greatest integer not
exceeding $(2\gamma+d)({1/p}-1).$ Then for any $f\in H^{p}(\mathbb{R}^{d})$
the Dunkl transform of $f$ satisfies the following Hardy's type
inequality
\begin{equation}\label{3.3.13}
\int_{\mathbb{R}^{d}}\frac{|\mathcal{F}_{\mathcal{D}}(f)(y)|^{p}}{|y|^{\sigma}}
d\mu_{k}(y)\leq C\|f\|^{p}_{H^{p}(\mathbb{R}^{d})},
\end{equation}
provide that
\begin{equation}\label{3.3.14}
(2\gamma+d)(2-p)\leq\sigma<2\gamma+d+p(N+1)
\end{equation}
where $C$ is a constant does not depends on  $f.$
\end{theorem}
\begin{Remarks}
\begin{verse}
\end{verse}
\begin{enumerate}
\item Note that the collection of all real $\sigma$ satisfying the
condition (\ref{3.3.14}) is a nonempty set since
$2\gamma+d+p(N+1)-(2\gamma+d)(2-p)>0.$
\item For the critical case $\sigma_0=(2\gamma+d)(2-p)$ has been extensively studied in \cite{Soltani}.
\item It would be interesting to know if this is the best possible improved.
\end{enumerate}
\end{Remarks}

\begin{proof}

Let $f=\sum_{j=0}^{\infty}\beta_{j}a_{j}\in H^{p}(\mathbb{R}^{d}),$
being element of $H^{p}(\mathbb{R}^{d})$ where $a_{j}$ are atoms.
Since $0 < p \leq 1 $ it follows
$$\int_{\mathbb{R}^{d}}\frac{|\mathcal{F}_{\mathcal{D}}(f)(y)|^{p}}{|y|^{\sigma}}
d\mu_{k}(y)\leq C
\sum_{j=0}^{\infty}|\beta_{j}|^{p}\int_{\mathbb{R}^{d}}\frac{|\mathcal{F}_{\mathcal{D}}(a_{j})(y)|^{p}}{|y|^{\sigma}}
d\mu_{k}(y).$$ In order to prove the theorem, it is enough to prove,
\begin{equation}\label{3.3.15}
\int_{\mathbb{R}^{d}}\frac{|\mathcal{F}_{\mathcal{D}}(f)(y)|^{p}}{|y|^{\sigma}}
d\mu_{k}(y)\leq C.
\end{equation}
Let us now take $\rho$ an arbitrary nonnegative real number, and
decomposing the left hand side of (\ref{3.3.15}) as

\begin{eqnarray*}
\int_{\mathbb{R}^{d}}\frac{|\mathcal{F}_{\mathcal{D}}(a_{j})(y)|^{p}}{|y|^{\sigma}}
d\mu_{k}(y) &=&\int_{|y|<\rho}\frac{|\mathcal{F}_{\mathcal{D}}(a_{j})(y)|^{p}}{|y|^{\sigma}}
d\mu_{k}(y)+ \int_{|y|\geq\rho}\frac{|\mathcal{F}_{\mathcal{D}}(a_{j})(y)|^{p}}{|y|^{\sigma}}
d\mu_{k}(y)\nonumber\\\nonumber\\
&:=& S_{1}+S_{2}.
\end{eqnarray*}
To estimate $S_{1};$ we may use Taylor's theorem in several
variables with integral's remainder for the function $ y\longmapsto E_{k}(ix,y),$ we obtain

$$E_{k}(ix,y)=\sum_{n=0}^{N}\frac{V_{k}(\prec iy,.\succ)^{n}(x)}{n!}+R_{N+1}(x,y) ,$$
where $$R_{N+1}(x,y)=\frac{1}{(N+1)!}\int_{0}^{1}(1-t)^{N+1}\Bigg[\int_{\mathbb{R}^{d}}
\prec iy,z\succ^{N+1}e^{t\prec iy,z\succ}d\Gamma_{x}(z)\Bigg]dt, $$
and $V_{k}$ is the intertwining operator (see, \cite{Dunkl1, Rosler1}), defined on $\mathbb{C}[\mathbb{R}^{d}]$
(the algebra of polynomial functions on $\mathbb{R}^{d}$) by
$$V_{k}(p)=\int_{\mathbb{R}^{d}}f(z)d\Gamma_{x}(z),~~~ x\in\mathbb{R}^{d}.$$
Since $\int_{\mathbb{R}^{d}}a_{j}(y)y^{\ell}d\mu_{k}(y)=0,$ for every $|\ell|\leq N $
where  $N=\Big[(2\gamma+d)(\frac{1}{p}-1)\Big],$ we can write
$$\mathcal{F}_{\mathcal{D}}(a_{j})(x)=\int_{B(0,r)}\Bigg[E_{k}(-ix,y)-
\sum_{n=0}^{N}\frac{V_{k}(\prec iy,.\succ)^{n}(x)}{n!}\Bigg]a_{j}(y)d\mu_{k}(y),~~~x\in\mathbb{R}^{d}.$$
Hence, from (\ref{3}), it follows that
$$|\mathcal{F}_{\mathcal{D}}(a_{j})(x)|\leq c_{k}\int_{B(0,r)}|R_{N+1}(x,y)||a_{j}(y)|d\mu_{k}(y).$$
But it is clear that
$$|R_{N+1}(x,y)|\leq \frac{1}{(N+1)!}[|x|.|y|]^{N+1} .$$
Now with the help of properties $(i), (ii)$ and  $(iii)$  for
$a(p,\infty,s)$-atoms of $H^{p}(\mathbb{R}^{d}),$ we get the
following results
\begin{eqnarray*}
|\mathcal{F}_{\mathcal{D}}(a_{j})(x)|&\leq& C \int_{B(0,r)}|x|^{N+1}|y|^{N+1}
\mu_{k}(B(0,r))^{-\frac{1}{p}}~d\mu_{k}(y)\\
&\leq& C r^{N+1+(2\gamma+d)(1-\frac{1}{p})}|x|^{N+1},
\end{eqnarray*}
where we have used (see, \cite{Xu}),
$$\mu_{k}(B(0,r))=\frac{r^{2\gamma+d}}{c_{k}2^{\gamma+d/2}\Gamma(\gamma+d/2+1)}.$$
Integrating with respect to the measure $d\mu_{k}$ over the domain
$0< |y|< \rho,$ we obtain
\begin{eqnarray*}
S_1:=\int_{|y|<\rho}\frac{|\mathcal{F}_{\mathcal{D}}(a_{j})(y)|^{p}}
{|y|^{\sigma}}d\mu_{k}(y) &\leq& C
r^{p(N+2\gamma+d+1)-(2\gamma+d)}\int_{|y|<\rho}|y|
^{p(N+1)-\sigma}~d\mu_{k}(y)\\
&\leq& C r^{-(2\gamma+d)+p(N+2\gamma+d+1)}\rho^{2\gamma+d+p(N+1)-\sigma}
\end{eqnarray*}
that is
\begin{equation}\label{3.3.16}
 S_{1}\leq Cr^{-(2\gamma+d)+p(N+2\gamma+d+1)}\rho^{2\gamma+d+p(N+1)-\sigma}
\end{equation}
provide that $\sigma<2\gamma+d+p(N+1)$ which follows from the inequality
(\ref{3.3.14}).

Now to estimate $S_{2},$ we may apply H\"{o}lder's inequality for
$q=\frac{2}{p}$ and Plancherel formula to get
\begin{eqnarray*}
S_{2} &\leq&
\Bigg(\int_{\mathbb{R}^{d}}\big(|a_{j}(y)|^{p}\big)^{\frac{2}{p}}
d\mu_{k}(y)\Bigg)^\frac{p}{2} \Bigg(\int_{|y|\geq\rho}|y|
^{\frac{2\sigma}{p-2}}d\mu_{k}(y)\Bigg)^{\frac{2-p}{2}}\\
&\leq & C \|a_{j}\|^{p}_{L^{2}({\mathbb{R}^{d}})}
\Bigg(\int_{y\geq\rho}|y|
^{\frac{2\sigma}{p-2}}d\mu_{k}(y)\Bigg)^{\frac{2-p}{2}}\\
&\leq &
C\|a_{j}\|^{p}_{L^{2}({\mathbb{R}^{d}})}\rho^{\frac{(2\gamma+d)(2-p)}{2}-\sigma},
\end{eqnarray*}
provide that $\frac{(2\gamma+d)(2-p)}{2}<\sigma,$  which is a consequence of
the left hand side of (\ref{3.3.14}).  Taking into account that
\begin{eqnarray*}
\|a_{j}\|_{L^{2}(\mathbb{R}^{d})}^{2}&=&\int_{\mathbb{R}^{d}}|a_{j}(y)|^{2}~d\mu_{k}(y)\nonumber\\
&\leq &\int_{B(0,r)}\big[\mu_{k}(B(0,r))\big]^{-\frac{2}{p}}d\mu_{k}(y)\nonumber\\
&\leq& C~ r^{-\frac{(2\gamma+d)(2-p)}{p}}.
\end{eqnarray*}
We obtain $\|a_{j}\|_{L^{2}(\mathbb{R}^{d})}^{p}\leq C~
r^{-\frac{(2\gamma+d)(2-p)}{2}}$ and hence,
\begin{equation}\label{3.3.19}
S_{2}\leq C~ r^{-\frac{(2\gamma+d)(2-p)}{2}}\rho^{\frac{(2\gamma+d)(2-p)}{2}-\sigma}.
\end{equation}

\underline{Case 1}. If $~~\sigma_0 =(2\gamma+d)(2-p).$ We put $\rho =
\frac{1}{r},~\forall~ r>0,$ then we have $S_{1} \leq C$ and
$S_{2}\leq C.$\\

\underline{Case 2}. If $~~(2\gamma+d)(2-p)<\sigma< (2\gamma+d)+p(N+1).$  We shall
discuss the cases $ 0<r<1$ and $
r\geq 1.$\\

Hence, in order to deal with the case $0<r<1,$ we need more precise
estimates, so we consider the set $\Upsilon_{\rho} ;$ the collection of all
numbers $\rho$ satisfying
\begin{equation}\label{3.3.17}
\frac{(2\gamma+d)(2-p)}{(2\gamma+d)(2-p)-2\sigma}\log(r)\leq
\log(\rho)\leq\frac{(2\gamma+d)-p(N+2\gamma+d+1)}{(2\gamma+d)+p(N+1)-\sigma}\log(r).
\end{equation}
To prove that the collection $\Upsilon_{\rho} $ above is nonempty set  it is
enough to prove that

\begin{equation}\label{C}
\frac{(2\gamma+d)(2-p)}{(2\gamma+d)(2-p)-2\sigma}\times
\frac{(2\gamma+d)+p(N+1)-\sigma}{(2\gamma+d)-p(N+1+(2\gamma+d))}\leq1
\end{equation}
which is a different formulation of the  hand side of
(\ref{3.3.14}), that is $(2\gamma+d)(2-p)\leq\sigma.$\\
Using the fact that $(2\gamma+d)+p(N+1)-\sigma>0 $ and the right hand side of
(\ref{3.3.17}) it follows that
\begin{equation}\label{3.3.18}
S_{1}\leq C  r^{-(2\gamma+d)+p(N+2\gamma+d+1)}\rho^{2\gamma+d+p(N+1)-\sigma}\leq C.
\end{equation}
Using the left hand side of (\ref{3.3.17}) and the fact that
$\frac{(2\gamma+d)(2-p)}{2}-\sigma <0,$ we obtain
\begin{equation}\label{3.3.20}
S_{2}\leq C.
\end{equation}
Combining (\ref{3.3.18}) and (\ref{3.3.20}) the result follows for
the case $0<r<1.$

Now, to deal with the case $ r \geq 1,$ we may take
\begin{equation}\label{3.3.21}
\rho =r^{\frac{2\gamma+d-p(N+1+2\gamma+d)}{2\gamma+d+p(N+1)-\sigma}}
\end{equation}
so, using the fact that $ r \geq 1,$ we obtain
\begin{equation}\label{3.3.22}
\rho =
 r^{\frac{(2\gamma+d)(2-p)}{(2\gamma+d)(2-p)-2\sigma}}.
\end{equation}
Combining (\ref{3.3.16}), (\ref{3.3.21}) and (\ref{3.3.19}) together
with (\ref{3.3.22}), the proof of the main theorem is completed.

\end{proof}


{\small
\begin{thebibliography}{}



\bibitem{CO} R. R. Coifman,
{\it A real-variable characterization of $H^p,$}  Studia Math. {\bf
51}, (1974), 269-274.

\bibitem{Jeu} M. F. E. De Jeu, {\it  The Dunkl transform,} Invent. Math. {\bf 113} (1993),  147–-162.

\bibitem{Dunkl} C. F. Dunkl, {\it Differential-difference operators associated to reflection groups,}
Trans. Amer. Math. Soc. {\bf 311} (1989),  167--183.

\bibitem{Dunkl1} C. F. Dunkl, {\it Integral kernels with reflection group invariance,} Canad. J. Math. {\bf 43} (1991), 1213--1227.

\bibitem{Dunkl2} C. F. Dunkl,  {\it Hankel transforms associated to finite reflection groups,}
 Contemp. Math. {\bf138} (1992),  123-–138.

\bibitem{FS} C. Fefferman and E. M. Stein,
{\it $H^p$ spaces of several variables,} Acta Math. {\bf 129},
(1972), 137-193.

\bibitem{S3} G. B. Folland and E. M. Stein, \emph{Hardy Spaces on Homogeneous
Groups,} Princeton University Press, Princeton, NJ, 1982.

\bibitem{GR} J. Garcia-Cuerva and J. Rubio de Francia,
\emph{Weighted Norm Inequalities and Related Topics,}  North
Holland, 1985.

\bibitem{Humphreys} J.E. Humphreys,  \emph{Reflection groups and Coxter groups}
Cambridge: Cambridge Univ. Press 1990, 1--204.
%\bibitem{kan}{Y. Kanjin,} {\it On Hardy-Type Inequalities and Hankel Transforms,}
%Monatshefte $f\ddot{u}$r Mathematik,  {\bf 127}, (1999), 311-319.

\bibitem{AM1} A. Rahmouni and M. Assal, {\it An improved Hardy's inequality associated with the
Euclidean Fourier transform}, Acta Mathematica Scientia, (2013), 1-6.

\bibitem{Rosler} M. R$\ddot{o}$sler, {\it  Positivity of Dunkl's intertwining operator,} Duke Math. J. {\bf 98} (1999), 445–-463.

\bibitem{Rosler1} M. R$\ddot{o}$sler, {\it A positive radial product formula for the Dunkl kernel,}
Trans. Amer. Math. Soc. {\bf355} (2003), 2413–-2438.

\bibitem{S2} E. M. Stein, \emph{Harmonic Analysis, real variable Methods, orthogonality and
oscillatory integrals,} Princeton Univ. Press, Princeton, NJ, 1993.


\bibitem{Soltani} F. Soltani,  \emph{Maximal Bochner–Riesz operators on Hardy-type spaces
in the Dunkl setting,} Integr. Transf. Spec. F., (2012), 1--15.

\bibitem{TW} M. H. Taibleson and G. Weiss,
{\it The molecular characterization of certain Hardy spaces,}
Ast\'erisque {\bf 77}, (1980), Soci\'et\'e Math. de France, Paris,
67-149.

\bibitem{Xu} S. Thangavelu and Y. Xu,  \emph{Convolution operator and maximal function for the Dunkl transform,} J. Anal. Math. {\bf 97} (2005), 25–-56.



\end{thebibliography}
}
\end{document}